\documentclass[a4paper,10pt]{article}
\title{  HOMOGENEOUS SPECTRUM,
DISJOINTNESS OF CONVOLUTIONS,
AND MIXING PROPERTIES OF DYNAMICAL SYSTEMS }
\author{V.V. Ryzhikov}
\begin{document}



\large
\vskip 0.5cm

\centerline{ HOMOGENEOUS SPECTRUM,
DISJOINTNESS OF CONVOLUTIONS,}
\centerline{AND MIXING PROPERTIES OF DYNAMICAL SYSTEMS\footnote{This work has been  published in the pilot issue of 
Selected Russian Mathematics, 
Vol.~1 (1999), no.~1, 13--24
 ( in fact not  available for readers).  }}
\centerline{{V.V.~Ryzhikov}}

     \vskip 0.5cm

\centerline{Abstract}
        \vskip .1cm

     \begin{quote}
{ In connection with Rokhlin's question on an automorphism
with a homogeneous nonsimple spectrum, we indicate a class of
measure-preserving maps $T$ such that
$T\times T$ has a homogeneous spectrum of multiplicity 2.
The automorphisms in question satisfy the condition
$\sigma\ast\sigma\perp \sigma$, where
$\sigma$ is the spectral measure of $T$.
We also show that there is a mixing automorphism
possessing the above  properties and their higher order analogs.
}\end{quote}

\vskip 20pt

\setcounter{section}{0}

\section{Introduction}

Let $T$ be an automorphism of a Lebesgue space
 $(X,  \mu)$, $\mu(X)=1$.
This automorphism induces a unitary operator
${\widehat{T}}\colon\, L_2(\mu)\to L_2(\mu)$, ${\widehat{T}} f(x)=f(Tx)$, and one can speak about
spectral measure of the operator ${\widehat{T}}$ and the multiplicity function
of the spectrum.

Rokhlin posed the question on  the existence of an automorphism with
a nonsimple homogeneous spectrum of finite multiplicity.
Automorphisms having nonsimple spectra
of finite multiplicity were constructed by Oseledets \cite{Oseledets}.
Katok \cite{Katok1} obtained  the following result:
for a generic collection of
maps $T$, the essential range of the multiplicity function of
$T\times T$ is either ${\mathcal M}_{T\times T}=\{ 2\}$  or
${\mathcal M}_{T\times T}=\{ 2, 4 \}$.
(Goodson and Lemanczyk noted \cite{GoodsonL} that
for $T\times T$ the spectrum multiplicity function does not take odd values.)
Katok conjectured that, for a generic $T$,
the automorphism $T\times T$ has a
homogeneous spectrum of multiplicity~2. This conjecture was confirmed
by Ageev and the author, see \cite{A},~\cite{Ry}.

 In this paper we study
automorphisms $T$ for which   ${\mathcal M}_{T\times T}=\{ 2\}$, i.e.
$T\times T$ has a homogeneous spectrum of multiplicity 2.
Namely, we consider automorphisms $T$ of simple spectrum such that
for $a\in (0,1)$,  the operator
$(aI+(1-a){{\widehat{T}}})$  (or the operator
$(1-a)(I+a{{\widehat{T}}}+a^2{{\widehat{T}}}^2+\dots)$)  belongs to the weak closure of
the powers of ${\widehat{T}}$ (\S\S 2, 3).
Using certain staircase constructions \cite{Adams},
one can obtain similar spectral properties for mixing automorphisms (\S 4).
In \S 5 we discuss disjointness of higher order convolutions.

For definitions,  we refer the reader  to \cite{G}, \cite{Stepin}.
We shall denote by $\sigma=\sigma_T$ the maximal
spectral type of the unitary operator ${\widehat{T}}$ acting on the  space
 $H =\{f\in L_2(\mu)\,\colon \int f d\mu = 0\}$.
It is well known that the spectral type of ${\widehat{T}}\otimes{\widehat{T}}$
is the measure  $\sigma + \sigma\ast\sigma$.  Note also that the convolution
$\sigma\ast\sigma$ is the maximal spectral type of the restriction of
${\widehat{T}}\otimes{\widehat{T}}$ on the space   $H\otimes H$.

Throughout we use the fact that
the disjointness of  $\sigma$ and  $\sigma\ast\sigma$
is equivalent to the absence of nonzero  operators intertwining
${\widehat{T}} |H$ with   ${\widehat{T}}\otimes {\widehat{T}} | H\otimes H$.
We note that the conditions $$({\widehat{T}}\otimes {\widehat{T}})
J=J{\widehat{T}}, \ \ \ J\neq 0$$
imply the existence of a complex  measure $\lambda\neq 0$
such that
$$  \widehat{\lambda}(n)=\langle {\widehat{T}}^n f|J^\ast g\rangle_{L_2(\mu)}   =
\langle {\widehat{T}}^n\otimes {\widehat{T}}^n  Jf|g\rangle_{L_2(\mu\times\mu) },
$$
hence $\lambda$ is a common component of the measures
$\sigma$ and $\sigma\ast\sigma$.

\section{The case
$\ \ {{\widehat{T}}}^{k_i}\to  \bigl(aI+(1-a){{\widehat{T}}}\bigr)$.}

Goodson \cite{G} noted that,  given an automorphism $T$,
the simplicity of the spectrum of the map  $R\colon\, X\times X \to X\times X,$
$$R(x,y)=(y,Tx),$$
implies that
$T\times T$  has a homogeneous spectrum of multiplicity~2. We shall use
this in the proofs below.
\medskip
\\
{\bf Theorem  2.1.} \it  
Let  $T$ be an ergodic automorphism and
let the weak convergence
$$\ \ {{\widehat{T}}}^{k_i}\to  (aI+(1-a){{\widehat{T}}})$$
hold
for some sequence ${k_i}\to \infty$ and some $a\in (0,1)$.
Then
\begin{itemize}
\item[$(1)$]
the spectral measure $\sigma$ of ${\widehat{T}}$ and the convolution
   $\sigma\ast\sigma$ are disjoint$;$

\item[$(2)$] if,  in addition,  $T$ has a simple spectrum, then
   $(T\times T)$ has a homogeneous spectrum of multiplicity~$2$.
\end{itemize}

\medskip
\rm

{Proof.}
(1) Suppose that a bounded operator $J$ intertwines ${{\widehat{T}}}$ and
$({{\widehat{T}}}\otimes {{\widehat{T}}})$, i.e.,
$$
J{{\widehat{T}}}= ({\widehat{T}}\otimes {\widehat{T}})J.
$$
Setting $b=1-a$  we obtain
$$
J(aI + b{{\widehat{T}}})=
\bigl((aI + b{{\widehat{T}}})\otimes (aI + b{{\widehat{T}}})\bigr)J,
$$
$$
\bigl(a(I\otimes I) + b({{\widehat{T}}}\otimes{{\widehat{T}}})\bigr)J =
 \bigl(a^2(I\otimes I) + ab({{\widehat{T}}}\otimes I) +ab
(I\otimes {\widehat{T}}) +b^2({\widehat{T}}\otimes {\widehat{T}})\bigr)J,
$$
$$
J + ({\widehat{T}}\otimes {\widehat{T}})J =
 (I\otimes {\widehat{T}})J +   ({\widehat{T}}\otimes I)J. \eqno (1)
$$
This implies that, for all  $i,j$,  we have
$$
({\widehat{T}}^i\otimes {\widehat{T}}^j)J + ({\widehat{T}}^{i+1}
\otimes {\widehat{T}}^{j+1})J -
({\widehat{T}}^{i}\otimes {\widehat{T}}^{j+1})J -
({\widehat{T}}^{i+1}\otimes {\widehat{T}}^{j})J =0.
$$
Now we obtain
$$  \sum_{0\leq i,j <n}
({\widehat{T}}^i\otimes {\widehat{T}}^j)J + ({\widehat{T}}^{i+1}
\otimes {\widehat{T}}^{j+1})J -
 ({\widehat{T}}^{i}\otimes {\widehat{T}}^{j+1})J -
({\widehat{T}}^{i+1}\otimes {\widehat{T}}^{j})J =
$$
$$
 0 = J + ({\widehat{T}}^n\otimes {\widehat{T}}^n)J -
 (I\otimes {\widehat{T}}^n)J -  ({\widehat{T}}^n\otimes I)J.
\eqno (2)
$$
Since $T$ is weakly mixing (it is not hard to check that
an eigenfunction of operator ${\widehat{T}}$ must be a constant function),
we have ${\widehat{T}}^{n_i}\to\Theta$, where $\Theta$ is
the orthogonal projection onto the space of
the constant functions.
Hence,
from  (2) we conclude
$$
J + (\Theta\otimes\Theta) J = (I\otimes\Theta   ) J + (\Theta\otimes I)J.
$$
Thus, $\mathrm{Im}(J)\perp H\otimes H$, which is equivalent to
the assertion that the zero  operator is a unique operator intertwining
${\widehat{T}} | H$ and ${\widehat{T}}\otimes{\widehat{T}} |H\otimes H$.
It is a well-known fact that the latter is equivalent to
$$
\sigma\ast\sigma\perp\sigma.
$$

 (2)  Let us show that the automorphism  $R$  has
a simple spectrum.
 Let $f$ be a cyclic   vector for
the operator ${\widehat{T}}$ acting on the space
$$
H=\Bigl\{h\in L_2(\mu)\colon\, \int h\, d\mu =0\Bigr\}.
$$
Let us prove   that
$V_{0,0}=f\otimes f$ is a cyclic vector for the restriction of
the operator ${\widehat{R}}$ to the space  $H\otimes H$.
We have to show that all vectors $V_{m,n}= T^mf\otimes T^nf$
belong to the space  $L$, the closure of the linear span of
the set  $\{R^iV_{0,0}\colon\, i\in {\bf Z}\}$.
(Note that  ${\widehat{R}} V_{m,n}=V_{n,m+1}$.)
We  have
$$  \ldots  V_{0,1},  V_{0,0},  V_{1,1},\ldots \in L.$$
From the relation
$ {{\widehat{T}}}^{k_i}\to  (aI+b){{\widehat{T}}})$  we obtain
$$
[(aI+b{\widehat{T}})\otimes (aI+b{\widehat{T}} )]V_{0,0}
= a^2V_{0,0}+b^2V_{1,1}
+abV_{0,1}+abV_{1,0}\ \in L.
$$
Hence
$$
V_{1,0}\in L, \ \ \ V_{0,2}= {\widehat{R}} V_{1,0}\in L.
$$

We shall assume that  $V_{0,i}\in L$ for $i=0,1, \dots, p$ and prove
that $V_{0,p+1}\in L$.
It is a well-known fact that the weak closure of the powers ${\widehat{T}}^n$
is a semigroup, hence it contains
the operators  $(aI+b{\widehat{T}})^p$.
It follows that the vector
$U_p=[(aI+b{\widehat{T}})^p\otimes (aI+b{\widehat{T}} )^p]V_{0,0}$
belongs to the space~$L$.
We write
$$
U_p= (a^pb^pV_{p,0}+a^pb^pV_{0,p}) +
\sum\limits_{m,n\colon\, |m-n|< p} c_{m,n}V_{m,n},
$$
where all $V_{m,n}$  belong to $L$ as $|m-n|<p$.
So we get $(V_{p,0}+V_{0,p})\in L$.
This implies that
$$
V_{p,0}\in L,  \ \ \
V_{0,p+1}={\widehat{R}} V_{p,0}\in L.
$$
Thus we have proved     that,   for all $m,n$,
$V_{m,n}=R^{2m}V_{0,n-m}\in L$, i.e., $L=H\otimes H,$
the restriction of ${\widehat{R}}$ to $H\otimes H$, has a simple spectrum.
Note that
the restriction of ${\widehat{R}}$ to  $(1\otimes H)+ (H\otimes 1)$
has also a simple spectrum (${\bf 1}\otimes f$ is a cyclic vector). Since
the action of the operator ${\widehat{R}}^2$ on $(1\otimes H)+ (H\otimes 1)$ and
the action of
${\widehat{R}}^2$ on $(H\otimes H)$ are disjoint, the same is true for
${\widehat{R}}$ and
we obtain that $R$ has a simple spectrum.


\bf Remarks. \rm 
     {\rm
(i) Theorem 2.1  was proved in part independently by Ageev,
see~\cite{A}.
He also
gave a solution to Katok's conjecture concerning the sets of spectral
multiplicities of $T\times T\times\dots \times T$ for generic automorphisms
$T$.

(ii) Katok and Stepin pointed out that there is a three interval
exchange transformation $T$ (in fact
an automorphism of the half-circle,
induced by some rotation of the unit circle)
possessing  the property of the
$(n,n+1)$-type approximation (see \cite{Katok1},
\cite{GoodsonR} for the definition). This property implies,
for some sequence $h(i)\to\infty$ and any $p>0$, the weak convergence
$$   T^{ph(i)} \to \frac{1}{2}(I+T^p).  \eqno  (\ast)$$
        Thus, we have
interval exchange transformations satisfying the conditions of
Theorem 2.1.
In addition, we obtain $\kappa$-mixing property   for our $T$ for
$\kappa =\frac{1}{2}$, i.e., for some sequence ${k_i}\to\infty$, one
has
$$
    {\widehat{T}}^{k_i}\to  \frac{1}{2}(I+ \Theta),
\eqno     (\ast\ast)
$$
where $\Theta$ is the orthogonal projection onto the space of
the constant functions.
This confirms  the corresponding   \it Oseledets conjecture \rm from
\cite{Os} (see also \cite{G}).
To see that ($\ast$) implies ($\ast\ast$) we use the fact that,
for some sequence $p_j\to\infty$,  we have ${\widehat{T}}^{p_j}\to \Theta$,
since $T$ is weakly mixing.
}

\section{The case $\ \ {{\widehat{T}}}^{k_i}\to  (1-a)
(I + a{\widehat{T}} +a^2{\widehat{T}}^2+\dots).$}

In this section  we consider automorphisms  which can be close to
the class of mixing automorphisms: if $a$ is close to~$1$, then,
for ergodic $T$, the operator
$(1-a) (I + a{\widehat{T}} +a^2{\widehat{T}}^2+\dots)$
will be close to $\Theta$.
\medskip
\\
{\bf Theorem 3.1.} \it
Suppose that $T$ is an ergodic automorphism such that,
for some sequence ${k_i}\to \infty$ and some $a\in (0,1)$, one has
the weak convergence
$$ {{\widehat{T}}}^{k_i}\to \
{(1-a)}(I+a{\widehat{T}} +a^2{\widehat{T}}^2+\dots).$$
Then
\begin{itemize}
\item[$(1)$] the spectral measure $\sigma$ of ${\widehat{T}}$ and the convolution
   $\sigma\ast\sigma$ are disjoint$;$

\item[$(2)$] if   $T$ has a simple spectrum, then  the automorphism
   $(T\times T)$ has a homogeneous spectrum of multiplicity~$2$.
\end{itemize} \rm

Proof.
 (1) We denote
 $$
P= (1-a)(I-a{\widehat{T}})^{-1}= {(1-a)}(I+a{\widehat{T}} +
a^2{\widehat{T}}^2+\dots).
$$
Let an operator $J\colon\, H\to H\otimes H$ satisfy the condition
$$
J{{\widehat{T}}}= ({{\widehat{T}}}\otimes {{\widehat{T}}})J.
$$
Since ${\widehat{T}}^{k_i}\to P$, we have
$$
JP=  (P\otimes P) J,
$$
$$
J{(1-a)}(I+a{\widehat{T}} +a^2{\widehat{T}}^2+\dots) =(P\otimes P)J,
$$
$$
(1-a)[I\otimes I+ a({\widehat{T}}\otimes{\widehat{T}}) +
a^2({\widehat{T}}\otimes{\widehat{T}})^2+\dots]J
=(P\otimes P)J,
$$
$$
(1-a)(I\otimes I - a({\widehat{T}}\otimes{\widehat{T}}))^{-1} J =
(1-a)^2(I -a{\widehat{T}})^{-1}\otimes (I-a{\widehat{T}})^{-1}J.
$$
For the commuting operators $A=\bigl(I\otimes I -
a({\widehat{T}}\otimes{\widehat{T}})\bigr)$
and $B=(I -a{\widehat{T}})\otimes (I-a{\widehat{T}})$,
the equality $A^{-1}J=B^{-1}J$
implies $AJ=BJ$, hence we obtain
$$
{(1-a)}(I\otimes I - a({\widehat{T}}\otimes{\widehat{T}})) J =
(I-a{\widehat{T}})\otimes (I-a{\widehat{T}})J,
$$
$$
[ I\otimes I +  {\widehat{T}}\otimes{\widehat{T}}]J =
[ I\otimes{\widehat{T}} + {\widehat{T}}\otimes I]J.
$$
As it was proved above, the latter implies
$J=0$, which is equivalent to $\sigma\ast\sigma\perp\sigma$.

 (2) As in the proof of Theorem {2.1}, we show that
the space $H\otimes H$ is a cyclic space for the
operator ${\widehat{R}}$.  Let us use the following well-known fact:
if there is a sequence of cyclic spaces
$$ C_0 \subseteq C_1\subseteq C_2 \subseteq C_3 \subseteq \dots,  $$
whose union is dense in $H\otimes H$, then  $H\otimes H$ is a cyclic space.
Let $C_n$, $n=0,1,2,\dots,$ be
the  cyclic space (for the operator ${\widehat{R}})$ generated by
the vector $W_n$, where
$$  W_n= [(I-a{\widehat{T}})^n \otimes (I-a{\widehat{T}})^n] f\otimes f
$$
($f$ is a cyclic vector for ${\widehat{T}}$ on $H$).

Since $C_n$ is invariant with respect to
${\widehat{R}}^2={\widehat{T}}\otimes {\widehat{T}}$,
from ${\widehat{T}}^{k(i)}\to P$  we obtain
$(P\otimes P)C_{n}\subseteq C_{n}$.
But  $(P\otimes P)W_n  =W_{n-1}$, thus,
$C_{n-1}\subseteq C_n$.

Let  $L$  be the closure of the union of the $C_n$'s.
We have to show that all vectors
$V_{n,0}$ belong to $L$ (this implies that all
 $V_{m,n}\in L$).

Let us check that $V_{1,0}\in L$.
We have
$$V_{0,1}, V_{0,0}, V_{1,1},\ W_1=
[(I -a{\widehat{T}} )\otimes (I -a{\widehat{T}})]V_{0,0}\in C_1.$$
Since $W_1= V_{0,0} + a^2V_{1,1} -a V_{0,1} - aV_{1,0}$,
we get  $V_{1,0}\in C_1$.

We obtain
${\widehat{R}} V_{1,0} = V_{0,2}\in C_1$.
Let us prove that  $V_{2,0}\in C_2$. The vector  $W_2$ is
a linear combination of the vectors $V_{i,j}$, $0\leq i,j \leq 2$.
We know that all these vectors except  $V_{2,0}$ are in
$C_1$. But $W_2\in C_2$ and $C_1\subseteq C_2$, hence
$V_{2,0}\in C_2$.

By induction, we obtain
$$
V_{n,0}\in C_n\subseteq L,\qquad
\forall\, n=0,1,2,\dots
$$
As in the proof of Theorem {2.1}
 we conclude that
$L=H\otimes H$ and obtain the simplicity of the spectrum of
the operator ${\widehat{R}}$.

Theorems {2.1} and  3.2 have been announced   in~\cite{Ry}.

\section{The  case of mixing $T$.}

A generalization of the above methods enables us to
prove the existence of a
mixing operator $T$ with the following properties:

\begin{itemize}
\item[1.] Spectrum of the symmetric product $T\odot T$ is simple:
 ${\mathcal M}_{T\odot T}=\{ 1\}$,
\item[2.] ${\mathcal M}_{T\times T}=\{ 2\}$,
\item[3.] $\sigma_T\ast\sigma_T\perp\sigma_T$.
\end{itemize}
We recall that ${T\odot T}$ denotes the restriction of
${T\times T}$ to the factor of $S$-fixed subsets of $X\times X$,
where  the map $S$ is defined as $S(x,y)=(y,x)$.
Note also that Property~1  implies  Property~2  and
Property~3, which is readily seen.

Adams \cite{Adams} proved the property of mixing for a large
class of rank~$1$ staircase constructions.  For some special
automorphisms of the Adams class, we are able to prove Property~1.
This is a positive answer to the corresponding question of
J.-P. Thouvenot.

Let us recall the definition of a staircase construction. Let
an automorphism $T$ admit, for any $n$, a partition
$\xi_n$ of $X$ into sets
\\
  $  B^1_n,\ \    TB^1_n,\ \            \ \    \dots\ \  
     \ \    T^{h_n-1}B^1_n,\ \    $   
 \\
 $   B^2_n,\     TB^2_n,\             \    \dots \  
     \    T^{h_n-1}B^2_n,\    T^{h_n}B^2_n,\     $
 \\
  $  B^3_n,\     TB^3_n,\             \    \dots\  
     \    T^{h_n-1}B^3_n,\    T^{h_n}B^3_n,\ \  
    T^{h_n+1}B^3_n,\   $
 \\
$    \dots   \dots   \dots   \dots   \dots   \dots  
    \dots   \dots \quad    $
 \\
   $ B^{r_n}_n,\    TB^{r_n}_n,\     \    \dots\  
     \    T^{h_n+1}B^{r_n}_n,\  
    T^{h_n+2}B^{r_n}_n,\    \dots,   
              T^{h_n+r_n-2}B^{r_n}_n,\    Y_n  $

such that
\\
    $$ B^2_n   = T^{h_n}B^1_n, \\
     B^3_n   = T^{h_n+1}B^2_n, \\
     \dots \\
     B^{r_n}_n   = T^{h_n+r_n-1}B^{r_n-1}_n$$
for all $n$, and  the sequence
of the partitions $\xi_n$ tends to the partitions into singletons
($\xi_n\to\varepsilon$).
If, in addition, for all $n$ we have
$$B^1_{n-1}= B^1_{n}{\sqcup} B^2_{n} {\sqcup} \dots{\sqcup} B^{r_{n}}_{n}, $$
we say that such an automorphism $T$ is  a {\it staircase construction}.
One can see that this construction is defined by $h_1$ and
the sequence $\{r_n\}$.

It is easily seen  that  $h_{n} +1$ is the number of
atoms in the partition $\xi_{n-1}$. Note also that
the set $Y_{n-1}$ is the union of the sets
   $$  T^{h_n}B^2_n,\ \        \\
     T^{h_n}B^3_n,\ \    T^{h_n+1}B^3_n,\ \      \\
     \dots    \dots    \dots    \\
     T^{h_n+2}B^{r_n}_n,\ \    \dots,\ \  
     T^{h_n+r_n-2}B^{r_n}_n,\ \    Y_n.$$

The Adams theorem asserts that in the case where $r_n \to\infty$ and
$\frac{(r_n)^2}{h_n} \to 0$, the {\it corresponding staircase}
construction $T$ is  mixing.

We can prove that there is a sequence $r_n \to\infty$ such that
$\frac{(r_n)^2}{h_n} \to 0$, and the corresponding
operator $T$ possesses the
property    ${{\mathcal M}_{T\odot T}=\{ 1\}}$.
\\
\medskip
{\bf Theorem 4.1.} \it
There is a staircase mixing construction $T$ possessing the
property  ${{\mathcal M}_{T\odot T}=\{ 1\}}$. \rm

We introduce special classes $St.C.(p_j, h_j)$ of staircase constructions,
and prove that the property
 ${\mathcal M}_{T_j\odot T_j}=\{ 1\}$  holds
for any $T_j \in St.C.(p_j, h_j)$.
Then we find a sequence of corresponding $T_j$ which approximates
sufficiently well some
staircase construction $T$ satisfying the conditions of
the Adams theorem.
Due to  the approximation procedure, the property
${\mathcal M}_{T\odot T}=\{ 1\}$  will be preserved.

\it Definition. \rm
We say that $T$ is in the class $St.C.(p, h)$
if $T$ is a staircase construction with a sequence $r_n$ satisfying the
following conditions:
\begin{itemize}
\item[1.] $\liminf_{n\to\infty} r_n =p$;
\item[2.] $\forall q, \ p\leq q\leq h, \ \ \exists n_i\to \infty$
such that  $r_{n_i +1}\to \infty$ and
$\forall \ \ i \ \ r_{n_i} =q$.
\end{itemize}

It follows from Condition~2 that the  operators
$$
P_q=\frac{1}{q}(I+{\widehat{T}}+{\widehat{T}}^2+\dots + {\widehat{T}}^{q-2}
+\Theta)
$$
are in $WCl(T)$ as $p\leq q\leq h$. (In fact
the powers ${\widehat{T}}^{-h_{n_i}}$ converge weakly to
$P_q$ as $r_{n_i} =q$ for all~$i$.)
\medskip
\\
{\bf Lemma.} \it
  Let $T$ be in $St.C.(p, h+2)$, and let $3p<h+2$.
Suppose that $B,TB,\dots, T^hB$
are disjoint measurable sets, and  denote by
$C_{B\times B}$
the cyclic space generated by vector
$\chi_{B}\otimes \chi_{B}$  under the action of the operator
${\widehat{T}}\otimes {\widehat{T}}$.
Then  the functions
$$F_{i,j}=\chi_{T^iB}\otimes \chi_{T^jB}+
\chi_{T^jB}\otimes \chi_{T^iB}, \ \ \ 0\leq i,j\leq h,
$$
belong to the cyclic space $C_{B\times B}$. \rm

{Proof.}
Since
$$F_{q,0}\in C_{B\times B} \Longrightarrow  F_{q+i,i}\in C_{B\times B}
$$
we only show that $F_{q,0}\in C_{B\times B}$.
We have
$$      P_{q+2}\chi_{B}\otimes  P_{q+2} \chi_{B},  \ \
P_{q+1} \chi_{B}\otimes  P_{q+1}\chi_{B},  \ \
         P_{q}\chi_{B} \otimes  P_{q}\chi_{ B}\   \in \ C_{B\times B}.
$$
Let
$$
G_{m}= \frac{m^2}{(1+\mu(B))^2} P_{m} \chi_{B}\otimes  P_{m}\chi_{B}.
$$
One can check that
$$F_{q,0}=
  {\mathrm Const}\,
  [G_{q+2} - G_{q+1}   -  ({\widehat{T}}\otimes {\widehat{T}}) G_{q+1} +
   ({\widehat{T}}\otimes {\widehat{T}}) G_{q}].
$$
Thus, for all $q\geq p$,  we have $F_{q,0}\in C_{B\times B}$, hence,
all the functions $F_{i,j}$  are in  $C_{B\times B}$.

Now  we  prove that, for all
$q < p$,  we have $F_{q,0}\in C_{B\times B}$  too.
Let us show this for $q=1$.
Since
$$
F_{p,0}= [({\widehat{T}}^p\otimes I)+(I\otimes {\widehat{T}}^p)]
\chi_{B\times B}
$$
and   $C_{B\times B}$  contains   $F_{p+1,0}$, we
obtain that the cyclic space generated by  $F_{p,0}$ contains
$$[({\widehat{T}}^p\otimes I)+(I\otimes {\widehat{T}}^p)]F_{p+1,0} =
[({\widehat{T}}^p\otimes I)+(I\otimes {\widehat{T}}^p)]
[({\widehat{T}}^{p+1}\otimes I)+(I\otimes {\widehat{T}}^{p+1})]
\chi_{B\times B}.$$
The latter can be represented as the sum  $F_{2p+1,0} + F_{p+1, p}$, where
$F_{2p+1,0}\in C_{B\times B}$. Since this sum is also in
$C_{B\times B}$, we get $F_{p+1, p}\in C_{B\times B}$
and, hence $F_{1, 0}$
belongs to  $C_{B\times B}$.
\medskip
\\
\bf Theorem 4.2.  \it
Let $T$ belong to $St.C.(p,\infty)$.
Then   ${\mathcal M}_{T\odot T}=\{ 1\}$.
\rm

{Proof.}
We consider the sequence $\xi_n$ and the corresponding sequence
of the sets $B_n=B^1_n$
(see the above definition of a staircase construction).
Let us consider the sequence of the cyclic spaces $C_{B_n\times B_n}$
for the operator ${\widehat{T}}\otimes{\widehat{T}}$.
Each symmetric function $F(x,y)$ can be approximated by
linear combinations of  $F_{i,j}$, where the  functions $F_{i,j}$
depend on  $n$ and $0\leq i,j\leq h_n$.
From the above lemma we obtain that
the symmetric product $L_2(\mu)\odot L_2(\mu)$ can be approximated
by the cyclic spaces   $C_{B_n\times B_n}$.
The Katok--Oseledets--Stepin
approach yields that  $L_2(\mu)\odot L_2(\mu)$ is a cyclic space
as well. Thus ${\mathcal M}_{T\odot T}=\{ 1\}$.

How to construct a mixing operator $T$  with the  property
${\mathcal M}_{T\odot T}=\{ 1\}$? In fact we prove only the existence.
A non-constructive description is this: we consider a sequence
$\{r_n \}$    such that,
for all $k$, we have $r_{2k+1} =  2k+1$, and
the sequence $r_{2k}$ \it extremely slowly \rm tends to the infinity  as
$k\to\infty$.
Now we explain why such a sequence $\{r_n \}$
can give the desired staircase
construction.

We shall consider a sequence
of automorphisms $T_j$ such that  $T_j$ is of class $St.C(p_j,\infty)$
on $X_j= \{x : T_j(x)\neq x \}$, where $\mu(X_j)\to 1$ and $p_j\to\infty$.
Given an automorphism $T_j$ with the corresponding sequence
$r_n^{(j)}$,
we choose a sufficiently large
number $N_j$ and change  this sequence only for
$n>N_j$.  We obtain a new sequence $r_n^{(j+1)}$ and the corresponding
construction $T_{j+1}$ such that $T_j$ differs from  $T_{j+1}$ only
on the very small set $Y_{N_j}$ (here $Y_{N_j}$ corresponds to the
automorphism $T_{j+1}$). Our
operator $T$ will be a limit automorphism, in fact
a staircase construction with $r_n\to\infty$. We can choose
$r_n \leq n$, then we obtain $\frac{(r_n)^2}{h_n} \to 0$,
which guarantees the property of mixing.

 The sequence $\{T_{j}\}$ is organized
so that, for any symmetric  functions $F(x,y)$,  the distance  between
$F$ and the space $C_{B_j\times B_j}$
(here the cyclic space is considered for the operator
${\widehat{T}}_j\otimes{\widehat{T}}_j$)
tends to zero.
However, it is possible  to obtain the same property  for the sequence
of the spaces $C_{B_j\times B_j}'$,
the cyclic spaces of the operator ${\widehat{T}}\otimes{\widehat{T}}$.
Indeed,
let   for some  $F$  ($F$
will be taken from a fixed finite collection of symmetric functions)
we have
$$
   \biggl\|
     F - \sum_{k=-N}^{N}a_k U_j^k \chi_{B_j\times B_j} \biggr\|
   \ < \ \varepsilon,
$$
where $U_j  = {\widehat{T}}_j\otimes{\widehat{T}}_j$.
However, we can ensure that
$$
   \biggl\|
     F - \sum_{k=-N}^{N}a_k U^k \chi_{B_j\times B_j} \biggr\|
   \ < \ 2\varepsilon
$$
for   $U  = {\widehat{T}}\otimes{\widehat{T}}$, since the measure of the
set $\{x\colon\ T_j(x)\neq T(x)\}$ can be made as small as desired
by an appropriate choice of $T_{j+1}, T_{j+2},\dots$~.

Thus, it is possible to find a construction such that,
for any fixed countable set of symmetric
functions $F$ and, hence,
for all  symmetric  functions
(due to the separability of $L_2\otimes L_2$),
the distance  between $F$ and the space $C_{B_j\times B_j}'$  tends to 0
as $j\to\infty$.

Thus   we  conserve the property ${\mathcal M}_{T\odot T}=\{ 1\}$
for some mixing staircase construction $T$.
\\
     {\bf Remark.}
     {\rm
It is worth noting that for mixing $T$ the property
$${\mathcal M}_{T\odot T}=\{ 1\}$$  implies
mixing of all orders. This
follows from  Host's theorem: 
if $\sigma_T\ast\sigma_T\perp\sigma_T$,
then the mixing automorphism $T$ possesses the multiple mixing property
\cite{Host}.
     }

\section{ On higher order properties.}

Disjointness of $\sigma\ast\sigma$ and  $\sigma$ in the case
$2{{\widehat{T}}}^{k_i}\to (I+{{\widehat{T}}})$
was obtained also by M.~Lemanczyk and
generalized by  F.~Parreau to  $\sigma^{\ast\ n} \perp \sigma$.
The problem  ``$\sigma^{*n} \perp \sigma^{*m}$?'' as $m>n>1$
remains open.
In some cases, we can obtain a
bit more than  $\sigma^{*n} \perp \sigma$.
\\
\medskip
\bf Theorem  5.1. \it

If,  for an ergodic automorphism $T$, there is a sequence ${k_i}\to\infty$
such that
$${{\widehat{T}}}^{k_i}\to P=(aI+b{{\widehat{T}}}), \ \ 1>a>b>0,$$
then $\sigma\ast\sigma\perp\sigma\ast\sigma\ast\sigma$,
where $\sigma$ is the spectral measure of the automorphism~$T$.
\rm\medskip

     {Proof.}
Suppose that some operator $J\colon\, L_2\otimes L_2\to
L_2\otimes L_2\otimes L_2$   satisfies the intertwining condition
      $$ J({\widehat{T}}\otimes{\widehat{T}})
= ({\widehat{T}}\otimes{\widehat{T}}\otimes{\widehat{T}})J.$$
We have to show that $J=0$.
Since
$${{\widehat{T}}}^{-k_i +1}\to Q=(bI+a{{\widehat{T}}}),$$
we obtain
 $$  J[P\otimes P -  Q\otimes Q]  =
   [P\otimes P\otimes P -  Q\otimes Q\otimes Q]J, \\
   J (a^2-b^2)[I\otimes I -  {\widehat{T}}\otimes{\widehat{T}}]  =
   [P\otimes P\otimes P -  Q\otimes Q\otimes Q]J,$$
$$  J (a^2-b^2) =
[P\otimes P\otimes P -  Q\otimes Q\otimes Q +(a^2-b^2)
{\widehat{T}}\otimes{\widehat{T}}\otimes{\widehat{T}}]J,$$
$$  (a^2-b^2-a^3+b^3)J =
[\dots +   (a^2-b^2-a^3+b^3){\widehat{T}}\otimes
{\widehat{T}}\otimes {\widehat{T}}] J,$$
$$  J =      [(\dots) +  {\widehat{T}}\otimes {\widehat{T}}\otimes
{\widehat{T}}] J,$$
where  $(\dots)$ is a linear combination of the operators
${\widehat{T}}\otimes I\otimes I$, $\dots$,
$I\otimes {\widehat{T}}\otimes {\widehat{T}}$.
Let, for  $h\in H\otimes H$, we have
$$  Jh =      [(\dots) +  {\widehat{T}}\otimes
{\widehat{T}}\otimes {\widehat{T}}] Jh. $$
We rewrite this for the spectral representation of ${\widehat{T}}$ as follows:
$$  f(x,y,z) =   [(\dots) +  x y z] f(x,y,z), $$
where $f$ is the image of $Jh$ in the spectral representation, and
$x,y,z$ belong to the unit circle.
We can see that, for any fixed $x,y$, there is  a unique point
$z$ such that  
$$ 0\neq f(x,y,z) =   [(\dots) +  x y z] f(x,y,z). $$
Since the measure $\sigma$ is continuous,
one has
$\sigma\otimes\sigma\otimes\sigma\bigl({\rm{support}}(f)\bigr)=0$,
i.e., $f=0$ in the space
$L_2(\sigma\otimes\sigma\otimes\sigma)$.
Thus we obtain
$$  Jh=0,\ \ J=0,\ \
\sigma\ast\sigma\ \perp\ \sigma\ast\sigma\ast\sigma.$$

\medskip
The following assertion is a natural generalization of
Theorem 4.1.
\medskip
\\
{\bf Theorem  5.2.} \it
There is a mixing automorphism $T$
with the following properties{\rm:}
\begin{itemize}
 \item[$(1)$] ${\mathcal M}_{T^{\odot n}}=\{ 1\}$,
 \item[$(2)$] ${\mathcal M}_{T^{\times n}}=\{ n, n(n-1),\dots, n!\}$,
 \item[$(3)$] $\sigma^{\ast k}\perp\sigma^{\ast m}$ for all $k>m>0$.
\end{itemize} \rm

The proof of Theorem {5.2} will be published in a separate paper. \footnote{See  V.V.Ryzhikov, “Weak limits of powers, simple spectrum of symmetric products, and rank-one mixing constructions”, Sbornik: Mathematics (2007), 198(5):733-754 }
\\
     {\bf Remark.}
     {\rm
Ageev \cite{A} proved property (2)  for generic
(non-mixing) automorphisms. In~\cite{PR}, property (3) has been
established for the well-known Chacon automorphism. We conjecture that
the Chacon automorphism has properties (1) and (2) as well.
}
\medskip

Author is indebted to G. R. Goodson, V. I. Oseledets, and J.-P.~Thouvenot
for the  questions which stimulated this work.

\vskip 1cm

 Author's address:
\\ Chair TFFA, Dept. of Mech. and Math.,
\\ Moscow State University, Moscow 119899.

\end{document}